\documentclass[12pt]{amsart}
\usepackage{enumerate}
\usepackage{enumitem}
\usepackage[colorlinks=true, linkcolor=blue, urlcolor=blue, citecolor=blue, anchorcolor=blue, pdfborder={0 0 0}]{hyperref}
\usepackage{url}
\usepackage{graphicx,color,colortbl}
\usepackage[dvipsnames,table,xcdraw]{xcolor}
\usepackage{cite}
\usepackage{amsthm, amsmath, amssymb}
\usepackage{mathtools}
\usepackage[top=45truemm, bottom=45truemm, left=30truemm, right=30truemm]{geometry}
\usepackage{subcaption}
\usepackage{cancel}
\usepackage{float}
\usepackage{tabularx}
\usepackage{makecell}
\usepackage{array}
\usepackage{ragged2e}
\usepackage{csquotes}
\usepackage{booktabs}

\definecolor{LightGreen}{RGB}{208,254,184}
\definecolor{LightGray}{RGB}{242,242,242}

\newcolumntype{P}[1]{>{\RaggedRight\hspace{0pt}}p{#1}}
\newcolumntype{L}{>{\begin{math}}l<{\end{math}}}%
\newcolumntype{C}{>{\begin{math}}c<{\end{math}}}%
\newcolumntype{R}{>{\begin{math}}r<{\end{math}}}%

\title[On families of elliptic curves $E_{p,q}:y^2=x^3-pqx$]{On families of elliptic curves $E_{p,q}:y^2=x^3-pqx$ that intersect the same line $L_{a,b}:y=\frac{a}{b}x$ of rational slope}

\author[Eldar Sultanow]{\href{https://orcid.org/0000-0001-5257-2236}{\includegraphics[scale=0.06]{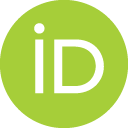}}\hspace{1mm}Eldar Sultanow}
\address{Eldar Sultanow\\Capgemini Deutschland GmbH\\Nuremberg, Germany}
\curraddr{}
\email{eldar.sultanow@capgemini.com}

\author[Anja Jeschke]{\href{https://orcid.org/0000-0001-7723-3986}{\includegraphics[scale=0.06]{orcid.png}}\hspace{1mm}Anja Jeschke}
\address{Anja Jeschke\\Capgemini Deutschland GmbH\\Hamburg, Germany}
\curraddr{}
\email{anja.jeschke@capgemini.com}

\author[A.\ Tfiha]{Amir Darwish Tfiha}
\address{Amir Darwish Tfiha\\Tishreen University\\Science Fuculty\\Syria}
\curraddr{}
\email{amirtfiha@tishreen.edu.sy}

\author[Madjid Tehrani]{\href{https://orcid.org/0000-0002-4838-5865}{\includegraphics[scale=0.06]{orcid.png}}\hspace{1mm}Madjid Tehrani}
\address{Madjid G. Tehrani\\George Washington University \\Washington, DC 20052, USA}
\curraddr{}
\email{madjid\_tehrani@gwu.edu}

\author[William J Buchanan]{\href{https://orcid.org/0000-0003-0809-3523}{\includegraphics[scale=0.06]{orcid.png}}\hspace{1mm}William J Buchanan}
\address{William J Buchanan\\Edinburgh Napier University\\Edinburgh, UK}
\curraddr{}
\email{b.buchanan@napier.ac.uk}

\subjclass[2010]{14H52}
\keywords{Elliptic Curves, Rational Points}

\begin{document}

\begingroup
\let\MakeUppercase\relax
\maketitle
\endgroup

\begin{abstract}
We investigate a special family of elliptic curves, namely $E_{p,q}:y^2=x^3-pqx$, where $p<q$ are odd primes. 
We study sufficient conditions for $p$ and $q$ so that the corresponding elliptic curve has non-trivial rational points. 
The number of sufficient conditions reduce to six. These six sufficient conditions relate to Polignac's conjecture, to the prime gap problem, the twin prime conjecture and to results from Green and Sawhney, and Friedlander and Iwaniec, e.g.
Additionally, we analyze the structures of the sufficient conditions for $p$ and $q$ by their graphical visualizations of the six sufficient conditions for $p,q \le 6997$. The graphical structures for the six sufficient conditions exhibit arc structures, quasi-linear arc segments, tile structures and sparsely populated structures.  
\end{abstract}

\section{Introduction}
\label{introduction}
The question of whether an elliptic curve has rational points or not has been occupying mathematicians for a while.  
The study~\cite{Daghigh_Didari_2015} proves that there is an upper bound of the rank of this family of elliptic curves. In addition, they prove that there are sufficient conditions on $p$ and $q$, for the elliptic curves to have rank two, three and four, and that there exist infinitely many such primes. Article~\cite{Burhanuddin_Huan_2014} shows that the additional restriction $p \equiv q \equiv 3 \mod{16}$ leads to the rank of the Mordell-Weil group $E_{pq}(\mathbf{Q})$ being less than or equal to 1, for example. 
Usual directions of studies take care about the rank of elliptic curves.
Our goal is to find structures of the conditions on the odd primes $p,q$ by visualization. We concentrate on sufficient conditions on the odd primes $p,q$ such that the corresponding elliptic curve has non-trivial rational points. 

\section{Candidate conditions}
\label{sec:candidate_conditions}
Let $L_{a,b}:y=\frac{a}{b}x$ be a linear function of rational slope where $a\in \mathbb{Z}, b\in \mathbb{N}$ and $(a,b)=1$. 
Let $E_{p,q}:y^2=x^3-pqx$ be an elliptic curve where $p<q$ are fixed odd primes. 
We aim to find sufficient conditions for $p$ and $q$ so that there are non-trivial rational points on the elliptic curve $E_{p,q}$. 
Inserting $L_{a,b}$ into $E_{p,q}$ leads to 
\begin{equation*}
x^3-\left(\frac{a}{b}\right)^2x^2-pqx=0.
\end{equation*}
This cubic function has two non-trivial solutions. They are given by
\begin{equation*}
x_{1,2}=\frac{1}{2}\left(\frac{a}{b}\right)^2\pm\frac{1}{2}\sqrt{\left(\frac{a}{b}\right)^4+4pq}. 
\end{equation*}
In order for $x_{1,2}$ to be rational, we need $a^4+4pqb^4$ to be a square: 
\begin{equation*}
a^4+4pqb^4=c^2
\end{equation*}
for some $c \in \mathbb{N}$. This equation can be factored as $4pqb^4=c^2-a^4=(c-a^2)(c+a^2)$. 
We want to find all possibilities of how to separate the factors of $4pqb^4$ into the two factors $c-a^2$ and $c+a^2$. \\
Counting the number of cases boils down to computing the number of divisors of $4pqb^4$. For $\tau(n)$ the divisor function with $n=p_1^{e_1}\cdots p_k^{e_k}$ and $p_k$ prime, we have the identity $\tau(n)=(e_1+1)\cdots(e_k+1)$, see~\cite[p.~34]{Remmert_Ullrich_2008}. If $b$ is prime, the number of cases would be
\begin{equation*}
\tau(2^2pqb^4)=(2+1)(1+1)(1+1)(4+1)=60.
\end{equation*}
If $b$ is not prime, the number of cases would be greater. 
We do not assume $b$ to be prime, but we will keep $b$ as a separate factor to keep the number of cases to consider at $60$.
Table~\ref{tab:cases_all} presents all 60 cases and their resulting conditions for $p$ and $q$. 
\begin{table}[htbp]
\centering
\begin{subtable}[b]{0.4\textwidth}
\centering
\begin{tabular}{l|lll}
 Case & $c-a^2$ & $c+a^2$ & Condition \\
 \toprule
  1 & $4pqb^4$ & 1 & $pq=\frac{1-2a^2}{4b^4}$ \\
  2 & $2pqb^4$ & 2 & $pq=\frac{1-a^2}{b^4}$ \\
  3 & $pqb^4$ & $4$ & $pq=\frac{4-2a^2}{b^4}$ \\
  4 & $b^4$ & $4pq$ & $pq=\frac{2a^2+b^4}{4}$ \\
  5 & $b^3$ & $4pqb$ & $pq=\frac{2a^2+b^3}{4b}$ \\
  6 & $b^2$ & $4pqb^2$ & $pq=\frac{2a^2+b^2}{4b^2}$ \\
  7 & $b$ & $4pqb^3$ & $pq=\frac{2a^2+b}{4b^3}$ \\
  8 & $2pb$ & $2qb^3$ & $p=qb^2-\frac{a^2}{b}$ \\
  9 & $4qb^4$ & $p$ & $p=2a^2+4qb^4$ \\
  10 & $2qb^4$ & $2p$ & $p=a^2+qb^4$ \\
  11 & $qb^4$ & $4p$ & $p=\frac{2a^2+qb^4}{4}$ \\
  12 & $qb^3$ & $4pb$ & $p=\frac{2a^2+qb^3}{4b}$ \\
  13 & $qb^2$ & $4pb^2$ & $p=\frac{2a^2+qb^2}{4b^2}$ \\
  14 & $qb$ & $4pb^3$ & $p=\frac{2a^2+qb}{4b^3}$ \\
  15 & $2qb$ & $2pb^3$ & $q=pb^2-\frac{a^2}{b}$ \\
  16 & $4pb^4$ & $q$ & $q=2a^2+4pb^4$ \\
  17 & $2pb^4$ & $2q$ & $q=a^2+pb^4$ \\
  18 & $pb^4$ & $4q$ & $q=\frac{2a^2+pb^4}{4}$ \\
  19 & $pb^3$ & $4qb$ & $q=\frac{2a^2+pb^3}{4b}$ \\
  20 & $pb^2$ & $4qb^2$ & $q=\frac{2a^2+pb^2}{4b^2}$ \\
  21 & $pb$ & $4qb^3$ & $q=\frac{2a^2+pb}{4b^3}$ \\
  22 & $2qb^2$ & $2pb^2$ & $p=q+\frac{a^2}{b^2}$ \\
  23 & $2b$ & $2pqb^3$ & $pq=\frac{a^2+b}{b^3}$ \\
  24 & $2b^2$ & $2pqb^2$ & $pq=\frac{a^2+b^2}{b^2}$ \\
  25 & $2b^3$ & $2pqb$ & $pq=\frac{a^2+b^3}{b}$ \\
  26 & $2b^4$ & $2pq$ & $pq=a^2+b^4$ \\
  27 & $4b$ & $pqb^3$ & $pq=\frac{2a^2+4b}{b^3}$ \\
  28 & $4b^2$ & $pqb^2$ & $pq=\frac{2a^2+4b^2}{b^2}$ \\
  29 & $4b^3$ & $pqb$ & $pq=\frac{2a^2+4b^3}{b}$ \\
  30 & $4b^4$ & $pq$ & $pq=2a^2+4b^4$ \\
\end{tabular}
\end{subtable}
\hfill
\begin{subtable}[b]{0.4\textwidth}
\centering
\begin{tabular}{l|lll}
 Case & $c-a^2$ & $c+a^2$ & Condition \\
 \toprule
  31 & 1 & $4pqb^4$ & $pq=\frac{2a^2+1}{4b^4}$ \\
  32 & 2 & $2pqb^4$ & $pq=\frac{a^2+1}{b^4}$ \\
  33 & 4 & $pqb^4$ & $pq=\frac{2a^2+4}{b^4}$ \\
  34 & $4pq$ & $b^4$ & $pq=\frac{b^4-2a^2}{4}$ \\
  35 & $4pqb$ & $b^3$ & $pq=\frac{b^3-2a^2}{4b}$ \\
  36 & $4pqb^2$ & $b^2$ & $pq=\frac{b^2-2a^2}{4b^2}$ \\
  37 & $4pqb^3$ & $b$ & $pq=\frac{b-2a^2}{4b^3}$ \\
  38 & $2qb^3$ & $2pb$ & $p=qb^2+\frac{a^2}{b}$ \\
  39 & $p$ & $4qb^4$ & $p=4qb^4-2a^2$ \\
  40 & $2p$ & $2qb^4$ & $p=qb^4-a^2$ \\
  41 & $4p$ & $qb^4$ & $p=\frac{qb^4-2a^2}{4}$ \\
  42 & $4pb$ & $qb^3$ & $p=\frac{qb^3-2a^2}{4b}$ \\
  43 & $4pb^2$ & $qb^2$ & $p=\frac{qb^2-2a^2}{4b^2}$ \\
  44 & $4pb^3$ & $qb$ & $p=\frac{qb-2a^2}{4b^3}$ \\
  45 & $2pb^3$ & $2qb$ & $q=pb^2+\frac{a^2}{b}$ \\
  46 & $q$ & $4pb^4$ & $q=4pb^4-2a^2$ \\
  47 & $2q$ & $2pb^4$ & $q=pb^4-a^2$ \\
  48 & $4q$ & $pb^4$ & $q=\frac{pb^4-2a^2}{4}$ \\
  49 & $4qb$ & $pb^3$ & $q=\frac{pb^3-2a^2}{4b}$ \\
  50 & $4qb^2$ & $pb^2$ & $q=\frac{pb^2-2a^2}{4b^2}$ \\
  51 & $4qb^3$ & $pb$ & $q=\frac{pb-2a^2}{4b^3}$ \\
  52 & $2qb^2$ & $2pb^2$ & $p=q-\frac{a^2}{b^2}$ \\
  53 & $2pqb^3$ & $2b$ & $pq=\frac{b-a^2}{b^3}$ \\
  54 & $2pqb^2$ & $2b^2$ & $pq=\frac{b^2-a^2}{b^2}$ \\
  55 & $2pqb$ & $2b^3$ & $pq=\frac{b^3-a^2}{b}$ \\
  56 & $2pq$ & $2b^4$ & $pq=b^4-a^2$ \\
  57 & $pqb^3$ & $4b$ & $pq=\frac{4b-2a^2}{b^3}$ \\
  58 & $pqb^2$ & $4b^2$ & $pq=\frac{4b^2-2a^2}{b^2}$ \\
  59 & $pqb$ & $4b^3$ & $pq=\frac{4b^3-2a^2}{b}$ \\
  60 & $pq$ & $4b^4$ & $pq=4b^4-2a^2$   
\end{tabular}
\end{subtable}
\caption{Conditions on $p$ and $q$ derived from 60 possible cases to split the term $4pqb^4$ into two factors $(c+a^2)$ and $(c-a^2)$.}
\label{tab:cases_all}
\end{table}

\section{Six sufficient conditions}
\label{sec:six_conditions}
\begin{figure}[H]
\includegraphics[clip, trim=0cm 17.3cm 0cm 0cm, width=1.00\textwidth, page=1]{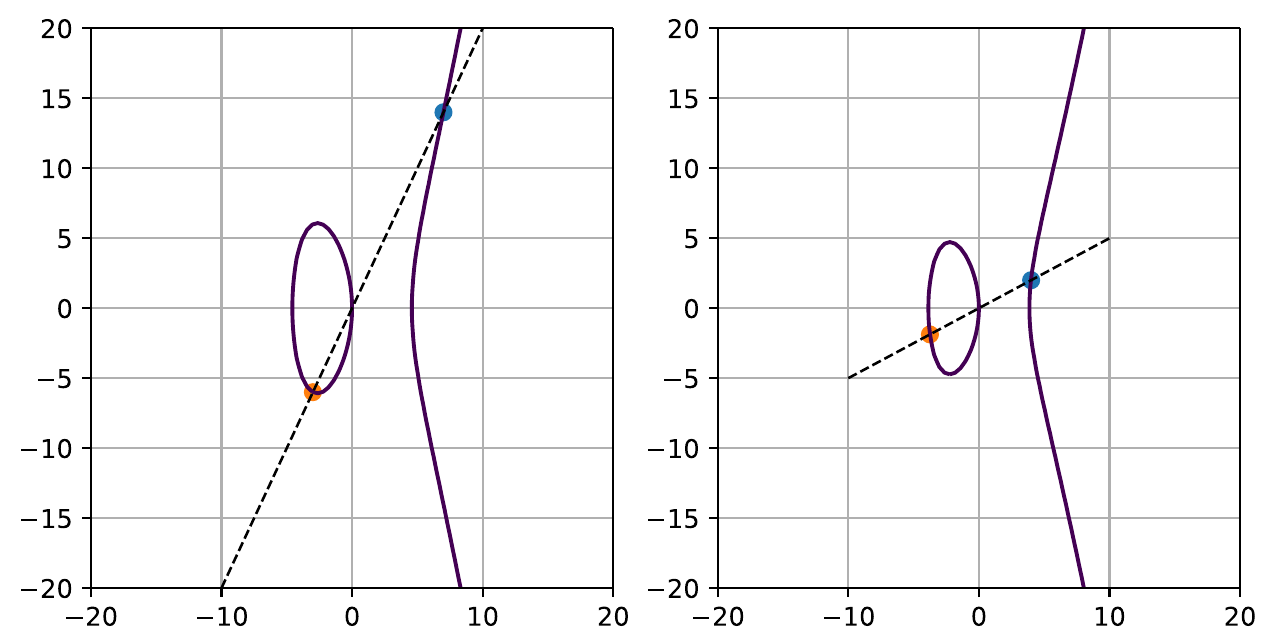}
\caption{Sample elliptic curves for case 40 (left) and case 56 (right) that are shown in Table~\ref{tab:six_conditions}. The two rational points are marked in orange and blue, respectively, for both cases. }
\label{fig:curves}
\end{figure}

The 60 candidate cases in Table~\ref{tab:cases_all} reduce to only six different cases that have solutions. The other 54 cases either have no solution (see Appendix~\ref{app:reasons_unsatisfiability}) or are redundant to one of the six cases (see Appendix~\ref{app:reasons_redundancy}). 
These six cases and their resulting corresponding conditions for $p$ and $q$ are listed in Table~\ref{tab:six_conditions}. Table~\ref{tab:six_conditions} also displays a sample combination of $a,b,p,q$, the resulting example elliptic curve and their rational points. As an illustration, Figure~\ref{fig:curves} shows the sample elliptic curves of case 40 and 56 together with their rational points.
\begin{table}[H]
\centering
\begin{tabular}{r|llll}
 Case & Condition & Samples & Resulting sample & Resulting rational points \\
  &  & a, b, p, q & elliptic curve & on sample elliptic curve \\
 \toprule
  17 & $q=a^2+pb^4$ & $8,3,3,307$ & $y^2=x^3-921x$ & $\left(\frac{307}{9},\frac{2456}{27}\right),(-27,-72)$ \\
  26 & $pq=a^2+b^4$ & $7,2,p,\frac{65}{p}$ & $y^2=x^3-65x$~\cite{LMFDB_135200.bq1} & $\left(\frac{65}{4},\frac{455}{8}\right),(-4,-14)$ \\
  32 & $pq=\frac{a^2+1}{b^4}$ & $1432,5,p,\frac{3281}{p}$ & $y^2=x^3-3281x$ & $\left(-\frac{1}{25},-\frac{1432}{125}\right),(82025,23491960)$ \\
  40 & $p=qb^4-a^2$ & $2,1,3,7$ & $y^2=x^3-21x$~\cite{LMFDB_14112.r1} & $(7,14),(-3,-6)$ \\
  47 & $q=pb^4-a^2$ & $5,2,3,23$ & $y^2=x^3-69x$~\cite{LMFDB_152352.ba1} & $(12,30),\left(-\frac{23}{4},-\frac{115}{8}\right)$ \\
  56 & $pq=b^4-a^2$ & $1,2,p,\frac{15}{p}$ & $y^2=x^3-15x$~\cite{LMFDB_14400.cq1} & $(4,2),\left(-\frac{15}{4},-\frac{15}{8}\right)$ \\
\end{tabular}
\caption{Sample elliptic curves including their LMFDB reference (if existing) and their rational points for all six cases that have solutions. }
\label{tab:six_conditions}
\end{table}

In the following, we give some additional details about the conditions of the six cases. We do not elaborate extensively because this is not the main focus of the paper. We give some links to known results in the literature, e.g. Polignac’s conjecture~\cite{dePolignac_1849}, the prime gap problem~\cite{Broughan_2021}, the twin prime conjecture~\cite{Cohen_Selfridge_1975, Murty_2013} and results from Green and Sawhney~\cite{GreenSawhney2024} and Friedlander and Iwaniec~\cite{Friedlander_Iwaniec_1997}.

\subsection{\texorpdfstring{Details of case 17: $\boldsymbol{q=a^2+pb^4}$}{Case 17}}
\label{sec:tab1_case17}
If $b$ is odd, the square $a^2$ is even, because $p$ and $q$ are odd. Hence, $a$ must be even. It is $(q-pb^4)\mod{4}=0$. The special case $b=1$ leads to $(q-p)\mod{4}=0$. 
According to Polignac's conjecture~\cite{dePolignac_1849} (which is not proved), there are infinitely many cases of two consecutive prime numbers with difference $n$ if $n$ is even~\cite[p. 295]{Bordelles_2020}.

\subsection{\texorpdfstring{Details of case 26: $\boldsymbol{pq=a^2+b^4}$}{Case 26}}
\label{sec:tab1_case26}
There are infinitely many primes of the form $p=a^2+b^4$, according to the theorem of Friedlander and Iwaniec ~\cite{Friedlander_Iwaniec_1997}. Friedlander and Iwaniec conjectured that the sum $a^2+4b^2$ also produces an infinite number of primes if $a$ and $b$ are prime numbers. Proving this is exactly what Green and Sawhney have achieved~\cite{GreenSawhney2024}. 
Up to a big integer $X$, there are $Y$ primes such that the product of two of them yields a semiprime of the form $a^2+b^4$. This observation aligns with Heath--Brown's general result on binary forms $G(x,y)$ of degree $d$, showing that the number of positive integers $n<X$ representable by such forms is of exact order $X^{2/d}$, provided $G(1,0)>0$ \cite[Theorem~8]{HeathBrown_2002}. 
For the specific case $G(x,y)=x^{2}+y^{4}$ (i.e., $d=4$), this implies a growth rate proportional to $X^{1/2}$, consistent with the density behavior established by Friedlander and Iwaniec~\cite{Friedlander_Iwaniec_1997}.

If $p\equiv3\pmod4$ or $q\equiv3\pmod4$, there is no solution~\cite[p.~21]{Aigner_Ziegler_2014}. Set $p=r^2+s^2$ and $q=u^2+v^2$. If $p\equiv1\pmod4$ and $q\equiv1\pmod4$ we obtain exact one solution $r>s>0$ for $p$ and one solution $u>v>0$ for $q$~\cite[p.~21]{Aigner_Ziegler_2014}. Hence, the product of both primes is $pq=(r^2+s^2)(u^2+v^2)=(ru+sv)^2+(rv-su)^2=(ru-sv)^2+(rv+su)^2$.
Consider $b^2=c$, other integer solutions for $pq=a^2+c^2=a^2+b^4$ do not exist, unless one of the four integers $ru+sv$, $|rv-su|$, $|ru-sv|$ and $rv+su$ is a perfect square.

\subsection{\texorpdfstring{Details of case 32: $\boldsymbol{pq=\frac{a^2+1}{b^4}}$}{Case 32}}
\label{sec:tab1_case32}
Every odd prime is of the form $4n+1$ or $4n+3$, i.e., every odd prime is congruent to $1$ or to $-1\pmod{4}$, see~\cite[p.~72]{Singh_2020},~\cite[p.~98]{LeVeque_1990}. The product of two odd prime numbers satisfies the following modular property:
\begin{itemize}
    \item If both primes satisfy $p\equiv1\pmod{4}$ and $q\equiv1\pmod{4}$, then their product is congruent with $p\cdot q\equiv 1\pmod{4}$.
    \item If one prime satisfies $p\equiv1\pmod{4}$ and the other $q\equiv3\pmod{4}$, then: $p\cdot q\equiv3\pmod{4}$.
    \item If both primes satisfy $p \equiv 3\pmod{4}$ and $q \equiv 3\pmod{4}$, then: $p \cdot q \equiv 1\pmod{4}$.
\end{itemize}
Thus, the product of two odd prime numbers is $1\pmod{4}$ or $3\pmod{4}$, depending on the residue classes of both prime numbers. The square of any even number is congruent to $0\pmod{4}$, while the squares of any odd number are congruent to $1\pmod{4}$, see~\cite[p.~21]{Aigner_Ziegler_2014}. Thus, quadratic residues modulo $4$ can only be $0$ or $1$. Considering the fourth power of $b$, we observe the following:
\[
b^4 \equiv 
\begin{cases}
0\pmod{4}, & \text{if } b \equiv 0, 2\pmod{4}, \\
1\pmod{4}, & \text{if } b \equiv 1, 3\pmod{4}.
\end{cases}
\]
Substituting this into the condition of case 32, we conclude: If $b^4 \equiv 0\pmod{4}$, then $a^2 \equiv 3\pmod{4}$, which is not possible. If $b^4 \equiv 1\pmod{4}$, then $a^2 \equiv 0\pmod{4}$, which means $a \equiv 0,2\pmod{4}$. This ensures that \( a \) and \( b \) are systematically constrained based on their modular properties.  

\subsection{\texorpdfstring{Details of case 40: $\boldsymbol{p=qb^4-a^2}$}{Case 40}}
\label{sec:tab1_case40}
If $b$ is odd, the square $a^2$ is even, because $p$ and $q$ are odd. Hence, $a$ must be even. It is $(qb^4-p)\mod{4}=0$. The special case $b=1$ results in $(q-p)\mod{4}=0$, which is the same as case 17 with $b=1$. Hence, according to Polignac's conjecture~\cite{dePolignac_1849} (which is not proved), there are infinitely many cases of two consecutive prime numbers with the difference n if n is even~\cite[p. 295]{Bordelles_2020}. If $b$ is even, it is $2^4|(p+a^2)$, which can be transformed to the problem of describing primes of the form $x^2+ny^2$. This was extensively elaborated by David A. Cox~\cite{Cox_2013}. 

\subsection{\texorpdfstring{Details of case 47: $\boldsymbol{q=pb^4-a^2}$}{Case 47}}
\label{sec:tab1_case47}
If $b$ is odd, the square $a^2$ is even, because $p$ and $q$ are odd. Hence, $a$ must be even. It is $(pb^4-q)\mod{4}=0$. The special case $b=1$ has no solution because it contradicts $q>p$.

\subsection{\texorpdfstring{Details of case 56: $\boldsymbol{pq=b^4-a^2}$}{Case 56}}
\label{sec:tab1_case56}
We can write $pq=(b^2-a)(b^2+a)$ which leads to $p=b^2-a$ and $q=b^2+a$ because $p<q$, assuming $a \in \mathbb{N}$ without loss of generality. Hence, $q-p=2a$. This is a special case of the prime gap problem. 
According to Polignac's conjecture~\cite{dePolignac_1849} (which is not proved), there are infinitely many cases of two consecutive prime numbers with the difference n if n is even~\cite[p. 295]{Bordelles_2020}. If $a=1$, this is well known as the twin prime conjecture~\cite{Cohen_Selfridge_1975, Murty_2013}.

\section{Graphical structures}
Using the six cases presented in Section~\ref{sec:six_conditions} we cover curves $E_{p,q}:y^2=x^3-pqx$, $p<q$ odd primes, by searching for tuples $(a,b)$, where $a\in \mathbb{Z}, b\in \mathbb{N}$ and $(a,b)=1$. 
In this search, we consider $p$, $q$ up to $6997$, and $b$ ranges from $1$ to $200$ and the first hit found for $a$ that meets the condition of the case considered is returned. The detailed code to visualize the structures is on GitHub~\cite{sultanow2024}. 

The pictures in this section should be read as follows: The quadratic structure is stored in matrix format. The entry in row $n$ and column $m$ is marked in yellow if there exists an $a\in \mathbb{Z}, b\in \mathbb{N}$, $(a,b)=1$ satisfying the condition of the case considered, e.g. $q=a^2+pb^4$, where $p$ is the $(n+1)$-th prime number and $q$ is the $(m+1)$-th prime number. The shift in the index occurs because we only consider odd primes for plotting. This means $n+1=\pi(p)$, $m+1=\pi(q)$. 
To improve readability, we have labeled the matrix axes with the actual prime numbers instead of their indices. Hence, the axes in all figures contain only odd primes. The plots show all odd primes p as a function of odd primes q, considering primes $p,q<6997$. 

The conditions of the cases are formulated in the trivial coordinates $(p,q) \in \mathbb{N}$. Hence, the descriptions of the cases 17, 40, 47 and 56 need to consider both the coordinates $(p,q) \in \mathbb{N}^2$ and the prime number coordinates $(\pi(p), \pi(q))$. 
The transformation factor between both coordinate systems is logarithmic due to the prime number theorem~\cite{Hadamard_1896, Poussin_1897} stating that the $i$-th prime number $p_i$ grows roughly like $p_i \sim i\log(i)$. Hence, $p \sim \pi(p) \log(\pi(p))$, $q \sim \pi(q) \log(\pi(q))$. The function $xlog(x)$ is locally curved, hence the figures of cases 17, 40, 47 and 56 show arcs instead of straight lines. The transformation factor between both coordinate systems is $\frac{\log(p_{\text{max}})}{\log(q_{\text{max}})}$, where $p_{\text{max}}, q_{\text{max}}$ are the maximal plotted values $p,q$ in each figure. 

The six cases reveal different structures. These structures are plotted and analyzed in more detail in the following subsections. 

\subsection{Case 17: Arc structure and quasi-linear arc segment}
\label{sec:structure_17}
\begin{figure}[H]
\includegraphics[clip, trim=0cm 0cm 0cm 0cm, width=1.00\textwidth]{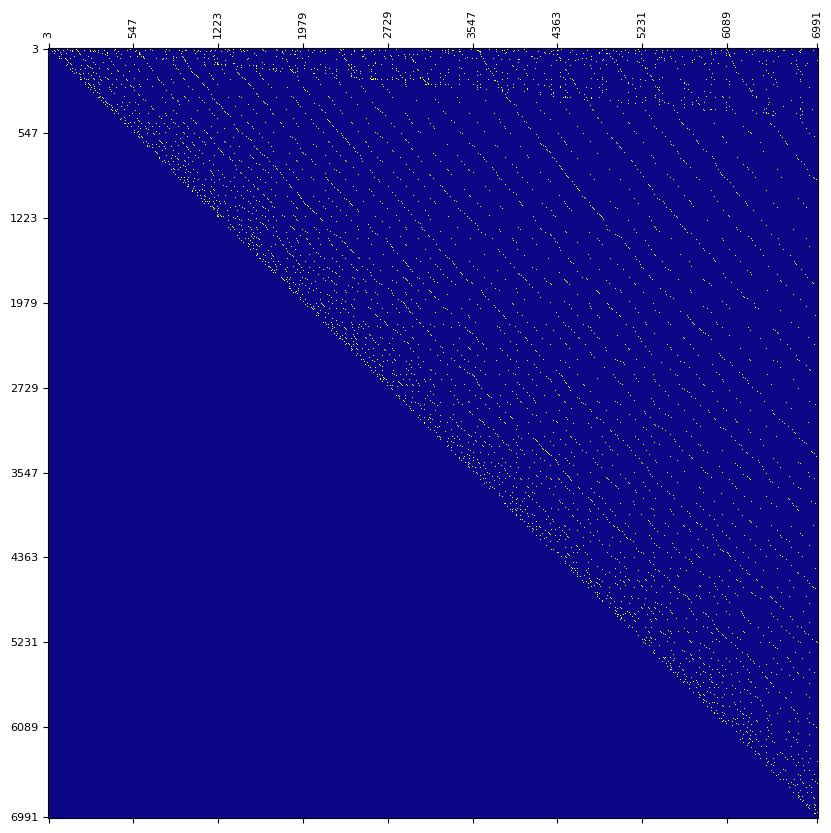}
\caption{Case 17 (Sections~\ref{sec:tab1_case17} and~\ref{sec:structure_17}) up to $p,q \le 6997$ showing arc structures. The plot shows odd primes p as a function of odd primes q. }
\label{fig:structure_case_17}
\end{figure}

Figure~\ref{fig:structure_case_17} displays the graphical structure of case 17 with its condition $q=a^2+pb^4$. For fixed $b$ and moderately small $|a|$, the solutions lie close to the affine lines $q \approx b^{4}p \quad(\text{up to the shift } a^{2})$ in the $(p,q)$-plane. These affine lines appear as arcs in Figure~\ref{fig:structure_case_17} due to the prime number coordinates used for plotting: Using the prime number theorem, the line $q \approx b^{4}p$ transforms approximately into $\pi(q)\log\pi(q) \;\approx\; b^{4}\,\pi(p)\log\pi(p)$, which is slightly curved in prime number coordinates. 

The parameter $b$ controls the ``slope'': $b=1$ yields arcs hugging the diagonal $q=p$, while larger $b$ push the arcs above that diagonal. The finite cut-off $p,q\le 6997$ truncates these arcs to visible segments.
Modular constraints described in Section~\ref{sec:tab1_case17} (e.g.\ $q - p b^{4} \equiv a^{2} \pmod 4$ with $a$ even when $b$ is odd) thin out these bands and create small lacunae along the arcs, exactly as seen in Figure~\ref{fig:structure_case_17}.

Additionally to these arc structures, a quasi-linear arc segment $q\approx 16p$ is visible for lower values of $p$. It is explained by setting $b = 2$ in the condition for moderately small $|a|$. 
A similar quasi-linear arc segment is visible in case 47, too. We refer to Section \ref{sec:structure_47} dealing with case 47 for a more detailed discussion. 

\subsection{Case 26: Tile structure}
\label{sec:structure_26}
\begin{figure}[H]
\includegraphics[clip, trim=0cm 0cm 0cm 0cm, width=1.00\textwidth]{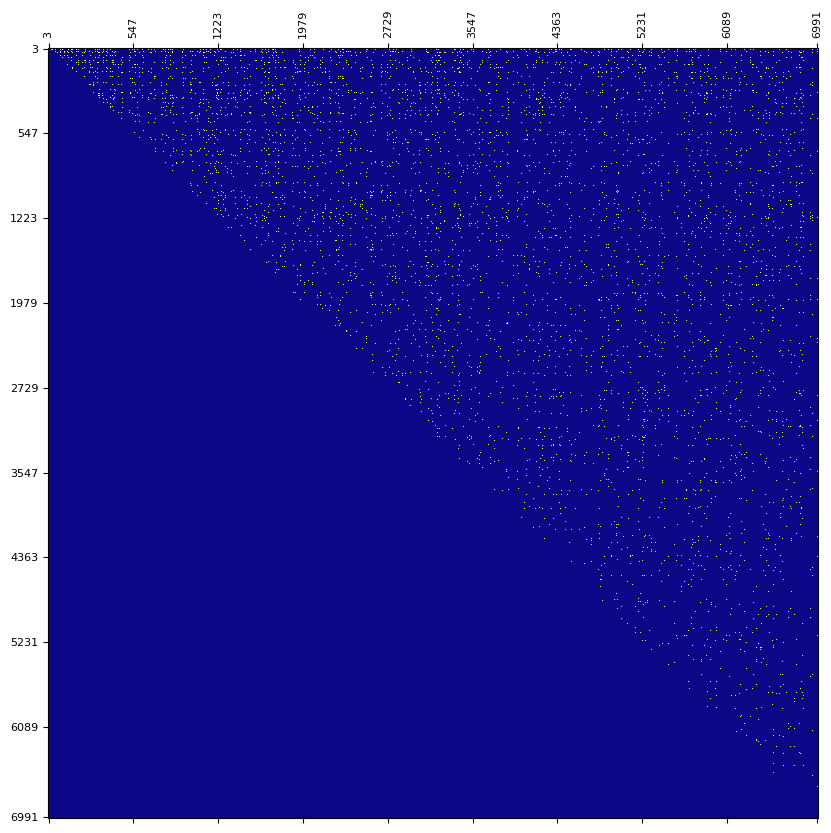}
\caption{Case 26 (Sections~\ref{sec:tab1_case26} and~\ref{sec:structure_26}) up to $p,q \le 6997$ showing tile structures. The plot shows odd primes p as a function of odd primes q.}
\label{fig:structure_case_26}
\end{figure}

Figure~\ref{fig:structure_case_26} displays the graphical structure of case 26 with its condition $pq=a^2+b^4$.
The tile structure shines through very slightly in Figure~\ref{fig:structure_case_26}. Such tiles may be due to the fact that for a certain prime number $p$ or $q$ the product $pq$ cannot be represented as a sum of $a^2+b^4$, or only with great difficulty. 
In detail, this phenomenon can be linked to deeper number-theoretic conjectures about representations of integers as sums of powers. In particular, the theorem from Friedlander and Iwaniec~\cite{Friedlander_Iwaniec_1997} shows that there are infinitely many primes of the form $p = a^2 + b^4$, but it does not guarantee that every product $pq$ of two primes can be expressed as $a^2 + b^4$ for integers $a$ and $b$. This structural limitation hints at the sparsity of solutions in the matrix, resulting in an almost checkerboard-like appearance. 
Furthermore, such expressions can also be viewed through the lens of the Hardy–Littlewood circle method~\cite{hardy1923partitio}, which predicts that representations of numbers as sums of squares and higher powers are often rare due to cancelation phenomena in exponential sums. 
Moreover, the Bunyakovsky conjecture~\cite{bunyakovsky1857} states that certain irreducible polynomials, such as $a^2 + b^4$, should produce primes infinitely often under certain conditions, but it does not cover composite products $pq$. Therefore, for many prime pairs $(p, q)$, there may simply be no integers $a$ and $b$ satisfying $pq = a^2 + b^4$, which explains the observed gaps and tile structures in the visualization.

\subsection{Case 32: Sparsely populated structure}
\label{sec:structure_32}
\begin{figure}[H]
\includegraphics[clip, trim=0cm 0cm 0cm 0cm, width=1.00\textwidth]{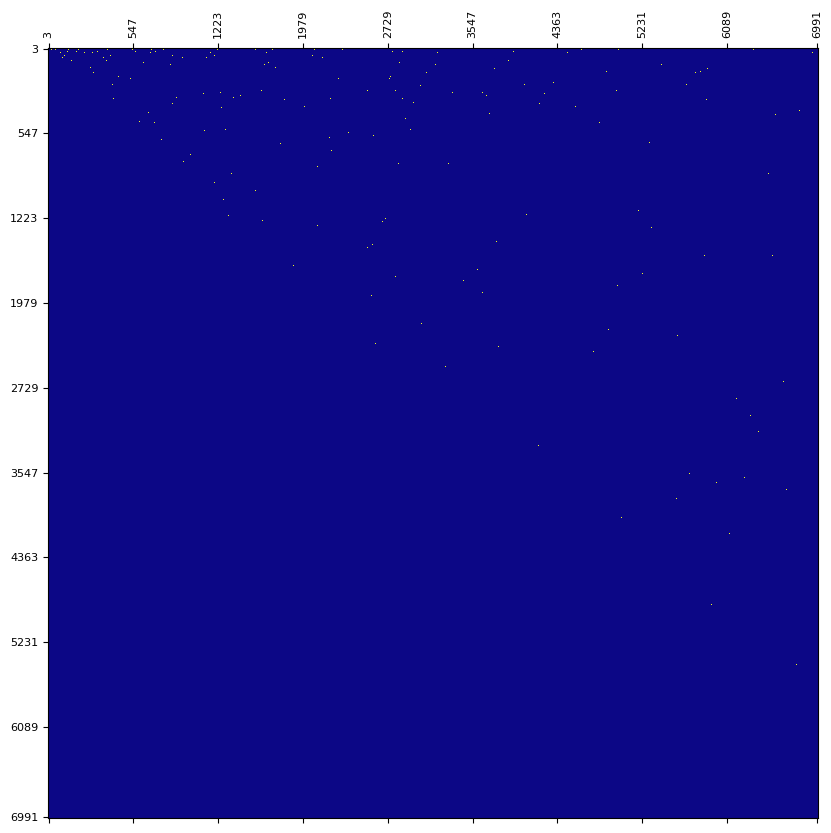}
\caption{Case 32 (Sections~\ref{sec:tab1_case32} and~\ref{sec:structure_32}) up to $p,q \le 6997$ showing sparsely populated solution area. The     plot shows odd primes p as a function of odd primes q.}
\label{fig:structure_case_32}
\end{figure}

Figure~\ref{fig:structure_case_32} displays the graphical structure of case 32. 
Solutions to the condition of case 32, namely $pq = \frac{a^2 + 1}{b^4}$, are rare: No solution has been discovered with all p, q, a, b being small. However, some small values of $p$ and $q$ appear in the graphical representation when at least one of the values of $a$ and $b$ is large. For example, a yellow dot in the image that is close to $(p,q)=(3,3)$ corresponds to $(p,q)=(5,337)$ with $(a,b)=(34522, 29)$. The smallest example for $b>1$ is $y^2=x^3-3281x$, which is the sample elliptic curve for case 32 shown in Table~\ref{tab:six_conditions}.
This rarity can be explained by the factorization properties of expressions of the form $x^2 + 1$. A well-known related case is the factorization of Fermat numbers, given by $F_n = 2^{2^n} + 1$. It is well known that every Fermat number $F_n$ is prime or has at least one very large prime factor~\cite{euler1732fermat, eulerArxiv}. This illustrates a general principle: Expressions of the form $x^2 + 1$ tend to have a very restricted set of prime factors or only a few, making it difficult for them to satisfy a product condition involving two small primes. 
Furthermore, the Bunyakovsky conjecture~\cite{bunyakovsky1857} on prime values of polynomials suggests that certain irreducible polynomials, such as $x^2 + 1$, take prime values only for a limited set of integers. Although there is no general proof of the conjecture, its implications suggest that for most values of $a$, the expression $a^2 + 1$ is unlikely to split into exactly two small prime factors, further restricting possible solutions to the given equation. 
These constraints, combined with the modular properties described in Section~\ref{sec:tab1_case32}, explain why valid pairs $(p,q)$ are extremely rare.

\subsection{Case 40: Arc structure}
\label{sec:structure_40}
\begin{figure}[H]
\includegraphics[clip, trim=0cm 0cm 0cm 0cm, width=1.00\textwidth]{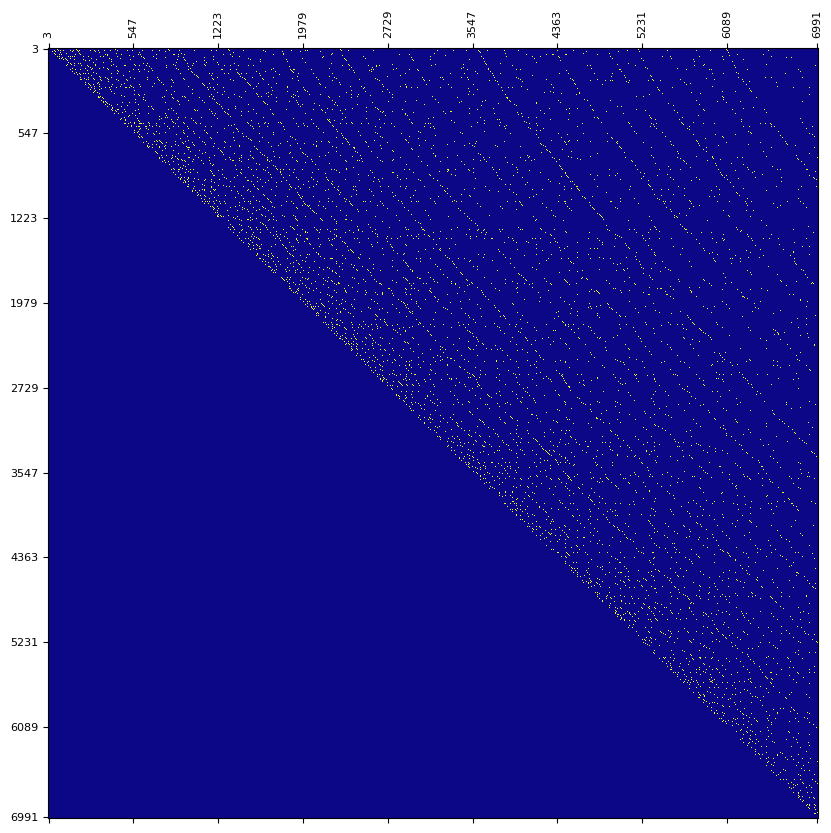}
\caption{Case 40 (Sections~\ref{sec:tab1_case40} and~\ref{sec:structure_40}) up to $p,q \le 6997$ showing arc structures. The plot shows odd primes p as a function of odd primes q.}
\label{fig:structure_case_40}
\end{figure}

Figure~\ref{fig:structure_case_40} displays the graphical structure of case 40 with its condition $p=qb^4-a^2$. The structure of cases 17 and 40 are similar, as the analysis in Sections \ref{sec:tab1_case17} and \ref{sec:tab1_case40} shows. 
For fixed $b$ and small $|a|$, valid pairs concentrate near the straight lines $p \approx b^{4} q \quad (\text{again up to the shift } -a^{2})$ in the $(p,q)$-plane. When plotting in prime number coordinates, these lines bend into arcs for the same reason as in case 17: Using the prime number theorem, the line $p \approx b^{4}q$ becomes $\pi(p)\log\pi(p) \;\approx\; b^{4}\,\pi(q)\log\pi(q)$, producing gently curved arc structures in Figure~\ref{fig:structure_case_40}.

The parity conditions summarized in Section~\ref{sec:tab1_case40} ($q b^{4} - p \equiv 0 \pmod 4$ and, for even $b$, the reduction to primes represented by $x^{2}+n y^{2}$) further sparsify these arcs, explaining the banded but porous appearance of the solution set. As $b$ increases, the main arcs shift away from the diagonal and become harder to distinguish within the plotting window $p,q\le 6997$, leaving only the lower-$b$ branches prominently visible.

\subsection{Case 47: Quasi-linear arc segment}
\label{sec:structure_47}
\begin{figure}[H]
\includegraphics[clip, trim=0cm 0cm 0cm 0cm, width=1.00\textwidth]{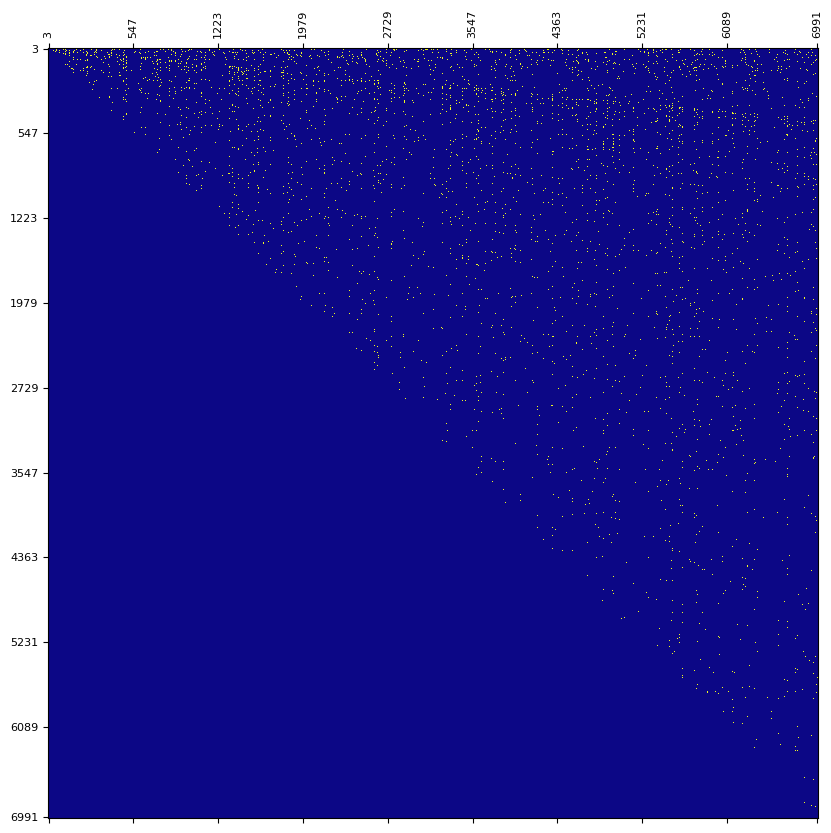}
\caption{Case 47 (Sections~\ref{sec:tab1_case47} and~\ref{sec:structure_47}) up to $p,q \le 6997$ showing quasi-linear arc segment structures. The plot shows odd primes $p$ as a function of odd primes $q$. A slight increase in the density of valid prime pairs can be observed for small $p < 437$.}
\label{fig:structure_case_47}
\end{figure}

\begin{figure}[H]
\includegraphics[clip, trim=0cm 0cm 0cm 0cm, width=1.00\textwidth]{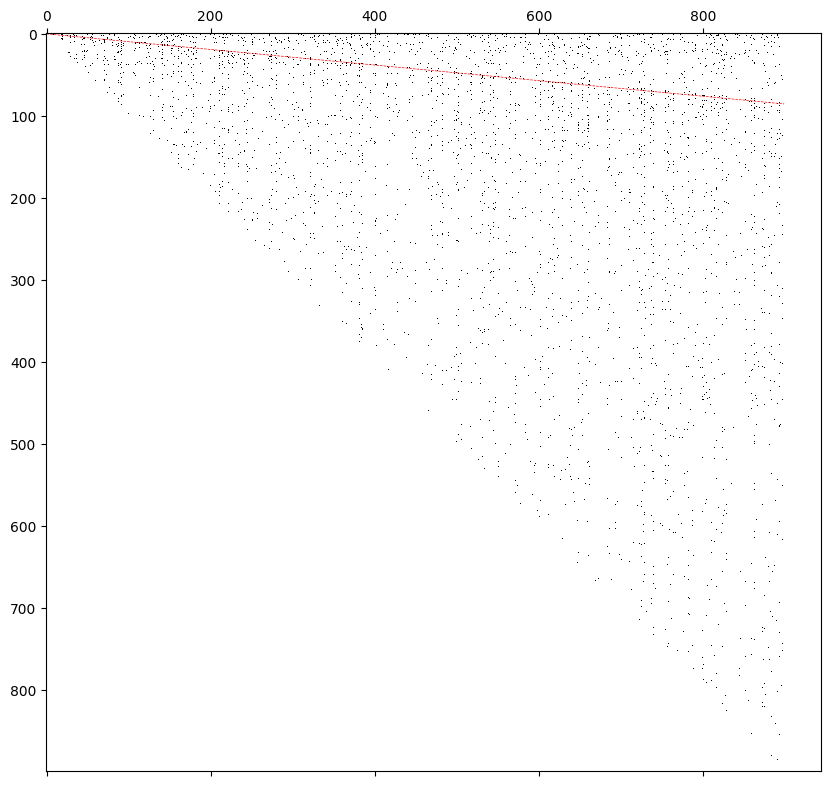}
\caption{Case 47: Pairs $(p,q)$ of odd primes $p<q$ solving $q=pb^4-a^2$ align along a line with slope $\approx 0.095$. The plot shows odd primes p as a function of odd primes q.}
\label{fig:structure_case_47_approx-line}
\end{figure}

Figure~\ref{fig:structure_case_47} displays the graphical structure of case 47 with its condition $q=pb^4-a^2$, reminiscent of a rotating pointer line on a submarine radar screen: We see a quasi-linear arc segment following the simple relation $q\approx16p$~\cite{5046353} for lower values of $p$. Figure~\ref{fig:structure_case_47_approx-line} shows the same visualization, but the quasi-linear arc segment highlighted in red. 
The quasi-linear arc segment is explained by setting $b = 2$ in the condition, which leads to $q=16p-a^2$. For small values of $|a|$, this simplifies to $q \approx 16p$, which means that valid prime pairs tend to accumulate on this quasi-linear arc segment or below this quasi-linear arc segment. Considering prime number coordinates, the slope of the quasi-linear arc segment in Figure~\ref{fig:structure_case_47} is approximately
\begin{align*}
\frac{\pi(p_{\text{max}})}{\pi(q_{\text{max}})} \approx \frac{\pi(\frac{1}{16}q_{\text{max}})}{\pi(q_{\text{max}})} =\frac{\pi(437)}{\pi(6997)} =\frac{84}{900} \approx 0.095.
\end{align*}
Hence, the fitted red line in the visualized data range up to $p_{\text{max}},q_{\text{max}}<6997$ in Figure~\ref{fig:structure_case_47_approx-line} is $\pi(p) \approx 10.5\pi(q)$ instead of $q\approx16p$. The transformation factor between the two coordinate systems is 
\begin{align*}
\frac{\log(p_{\text{max}})}{\log(q_{\text{max}})} \approx\frac{\log(\frac{1}{16}q_{\text{max}})}{\log(q_{\text{max}})} = \frac{\log(\frac{1}{16})+\log(q_{\text{max}})}{\log(q_{\text{max}})} \approx \frac{-1.2+3.84}{3.84} \approx 0.69 \approx \frac{10.5}{16}.
\end{align*}
When looking at larger numbers $(p_{\text{max}},q_{\text{max}})$, the transformation factor $\frac{\log(p_{\text{max}})}{\log(q_{\text{max}})}$ will approach $1$, and thus the quasi-linear arc segment will gradually curve toward the slope $\frac{1}{16}$~\cite{5046353}. 
Hence, quasi-linear arc segment in Figure~\ref{fig:structure_case_47} are a small extract of arc structures, similar to the observations w.r.t. cases 17 and 40. 

Additional quasi-linear arc segment structures should emerge for higher values of $b \in \mathbb{N}$. For $b=3$ and $b=4$, for example, we expect quasi-linear arc segments with $q\approx81p$ and $q\approx256p$, etc. These structures should be observable above the quasi-linear arc segment for $b=2$ by analyzing a sufficiently large dataset concerning $p$ and $q$. Our choice $p,q<6997$ makes these additional quasi-linear arc segment structures hard to distinguish from each other, although Figure~\ref{fig:structure_case_47} shows a slight increase in the density of valid prime pairs in the area $p<6997/16 \approx <437$. This becomes more apparent under logarithmic scaling or higher contrast suggesting a deeper regularity in how prime numbers align under this equation and motivating further investigation into their distributional properties.

\subsection{Case 56: Arc structure}
\label{sec:structure_56}
\begin{figure}[H]
\includegraphics[clip, trim=0cm 0cm 0cm 0cm, width=1.00\textwidth]{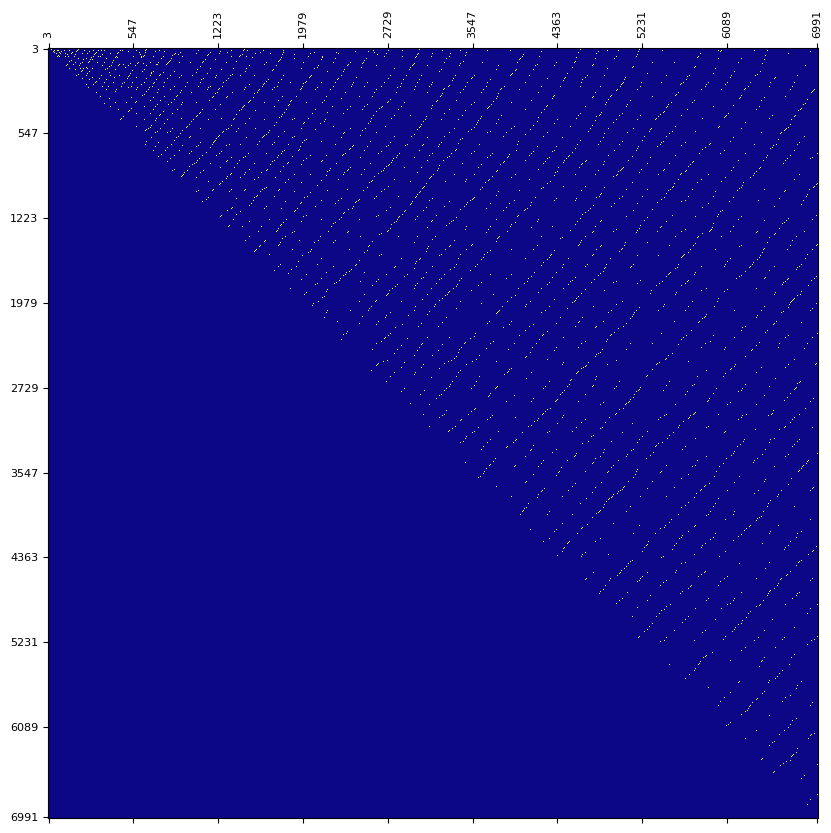}
\caption{Case 56 (Sections~\ref{sec:tab1_case56} and~\ref{sec:structure_56}) up to $p,q \le 6997$ showing arc structures. The plot shows odd primes p as a function of odd primes q.}
\label{fig:structure_case_56}
\end{figure}

\begin{figure}[H]
\includegraphics[clip, trim=0cm 0cm 0cm 0cm, width=1.00\textwidth]{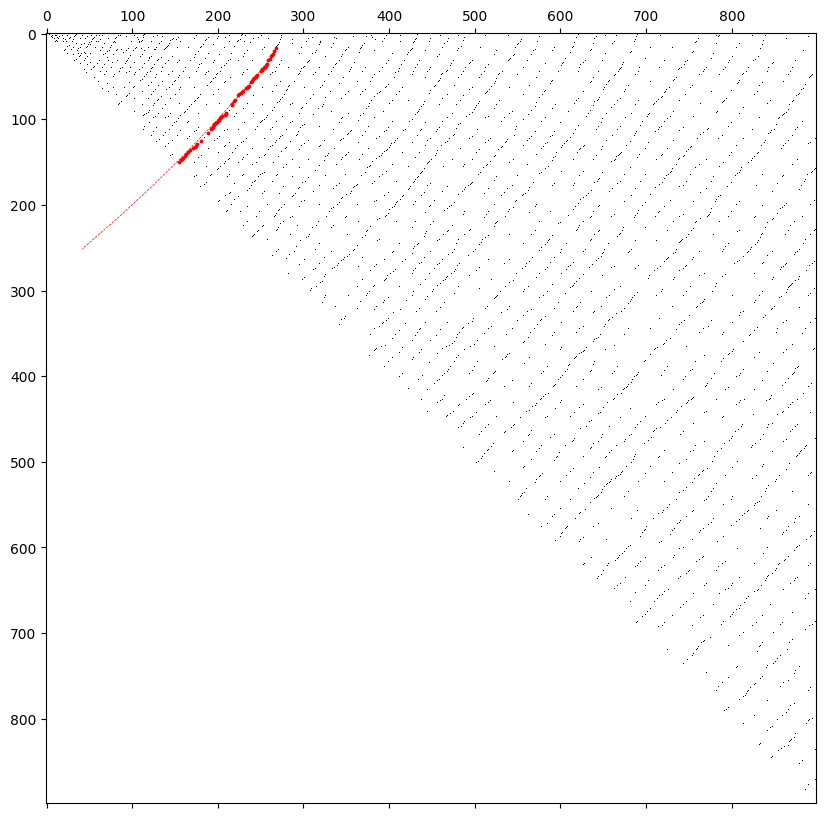}
\caption{Cases 56: Highlighted red points following a circular arc centered at $(-1051, -1051)$ with radius $R = 1700$. The plot shows odd primes p as a function of odd primes q.}
\label{fig:structure_case_56_circle-sample}
\end{figure}

Figure~\ref{fig:structure_case_56} displays the graphical structure of case 56 with its condition $pq=b^4-a^2$, showing slightly bent arcs. 
To investigate further, we extracted one of these arcs and found that a simple circle equation $(x - x_0)^2 + (y - y_0)^2 = R^2$ with a fixed center $x_0,y_0$ and radius $R$ captures almost all points of such an arc. 
For example, let us choose $x_0 = -1051, y_0 = -1051, R = 1700$. We almost perfectly capture a significant portion of these points as highlighted in red in Figure~\ref{fig:structure_case_56_circle-sample}. 

This visual phenomenon can be explained straightforwardly by taking a closer look~\cite{5035916} at condition $pq=b^4-a^2$ of case 56: This condition can be factored as $b^4-a^2=(b^2+a)(b^2-a)$, which suggests that for prime numbers $p$ and $q$, it is $q = b^2 + a$ and $p = b^2 - a$ due to $p < q$. By rearranging, we obtain:
\[
a = \frac{q - p}{2}, \quad b = \sqrt{\frac{p + q}{2}}.
\]
Thus, solutions exist when $\frac{p + q}{2}$ is a perfect square. These solutions lie along the lines $p+q=2b^2$ in the coordinates $(p,q) \in \mathbb{N}^2$. Plotting in prime number coordinates distorts the lines to arcs: Considering the prime number theorem, the lines $p+q=2b^2$ transform to $\pi(p)\log \pi(p)+\pi(q)\log \pi(q)=2b^2$, generating slightly curved arc structures in Figure~\ref{fig:structure_case_56}.

\newpage

\newpage
\section{Discussion and Outlook}
\label{outlook}
So far, we have inferred the conditions that two distinct odd primes $p,q$, $p<q$ must satisfy so that the elliptic curve $E_{p,q}:y^2=x^3-pqx$ can have rational points. There exist many curves that fulfill at least one of the six conditions and therefore have rational points.
Similar curves of the form $E_n: y^2=x^3-n^2x$ with $n\neq0$ a rational number are well known and exhibit the following properties: 
\begin{enumerate}
\item The rational number $n$ is a congruent number if and only if there is a point $P=(x,y)\in E_n(\mathbb{Q})$ with $y\neq0$, see~\cite[169]{Milne_2021},~\cite[144]{Stein_2009}. 
\item The only rational points of finite order in $E_n$ are four points of order $2$, namely $\mathcal{O}$ (the point at infinity), $(0,0)$ and $(\pm n,0)$, see~\cite[44]{Koblitz_1993}. In other words $E_n(\mathbb{Q})_{\text{tors}}=\{\mathcal{O},(0,0),(n,0),(-n,0)\}$. 
\end{enumerate}
Combining these two statements leads to the conclusion that the rational number $n$ is a congruent number if and only if the elliptic curve $E_n=y^2=x^3-n^2x$ has a positive rank. This is true, because the existence of a point $P=(x,y)\in E_n(\mathbb{Q})$ with $y\neq0$ is equivalent to $P$ is not in $E_n(\mathbb{Q})_{\text{tors}}$. Hence $P$ has infinite order, and in addition the elliptic curve has a positive rank. Researchers have also already addressed congruent numbers with many prime factors using the curve $y^2=x^3-n^2x$, see for example~\cite{Tian_2012}. 
Open questions w.r.t. the elliptic curve $E_{p,q}:y^2=x^3-pqx$ are, e.g.:
\begin{enumerate}
\item What must the criteria for $p$ and $q$ be for this family of elliptic curves to have a positive rank? 
\item Under which circumstances for given primes $p,q$ one of our six conditions is sure to yield an elliptic curve with a positive rank? 
\end{enumerate}
The criterion that $p<q$ are odd primes is not sufficient because two counterexamples are $p=5, q=7$, and $p=3, q=11$. The corresponding elliptic curves have rank zero, and their torsion generators are trivial, see LMFDB~\cite{LMFDB_78400.ga1, LMFDB_34848.ci1}. These questions may be answered in future investigations. 
Further interesting avenues for future research include:
\begin{itemize}
    \item Fixing a line and exploring families of curves where $pq$ is a congruent number. Visualizing the corresponding $(p,q)$ pairs may reveal underlying structures.
    \item Fixing a curve and investigating families of lines that intersect it at rational points, potentially uncovering deeper arithmetic properties.
    \item Analyzing the visualized structures in greater detail and interpreting them to draw meaningful conclusions about the distribution of prime numbers.
\end{itemize}

\section{Acknowledgements}
\label{acknowledgements}
We are grateful for the kind assistance and useful input from mathematics communities like the \hyperlink{https://math.stackexchange.com}{Stack Exchange Network} or the \hyperlink{https://www.matheboard.de}{MatheBoard Community}. As an example let us mention \hyperlink{https://math.stackexchange.com/users/232/qiaochu-yuan}{Qiaochu Yuan}, \hyperlink{https://math.stackexchange.com/users/30382/servaes}{Servaes} and \hyperlink{https://www.matheboard.de/profile.php?userid=39859}{HAL 9000}, for whose instant help in matters of number theory we are very grateful as for the help by many other pleasant people from the math communities.

We owe the graphical illustration of the structures to \href{https://stackoverflow.com/users/941531/arty}{Arty}. He implemented the C++ programs that generate the needed data, which can be checked out on GitHub~\cite{sultanow2023}.

\newpage
\appendix
\section{Unsolvable conditions}
\label{app:reasons_unsatisfiability}
We provide reasons for all conditions in Table~\ref{tab:cases_all} that cannot be solved: 
\begin{enumerate}[label=\alph*)]
\item Cases 1, 2, 3, 37, 53 and 54 are impossible to solve: This is due to the fact that $pq$ is an integer requires the fraction's numerator to be larger than the denominator. Case 3 requires the additional argument that $p,q$ need to be odd and hence, $a=b=p=1$ and $q=2$ do not provide a valid solution. 
\item Cases 4, 5, 6, 11, 12, 13, 18, 19, 20, 27, 28, 33, 34, 35, 36, 41, 42, 43, 48, 49, 50, 57, 58 are impossible to solve: The conditions of these cases require $b$ to be even and as a consequence $a$ to be odd (since $a\in\mathbb{Z}$, $b\in\mathbb{N}$, $(a,b)=1$). However, dividing the conditions by 2 leads to the conclusion that a must be even which is a contradiction. Let us take for example case 4 with the condition $4pq=2a^2+b^4$ that leads to $2pq=a^2+\frac{b^4}{2}$ when dividing it by $2$ and thus to the contradictory requirement $a$ is even. 
\item Case 7 is unsolvable: We know that $b$ must be even and by substituting $b$ with $2t$ we have to solve $16t^3pq=a^2+t$ which is $(16t^2pq-1)t=a^2$. Since $t$ is coprime with $16t^2pq-1$, we conclude that $16t^2pq-1$ is a perfect square, which is impossible by an argument $\bmod4$: Recall that if $x$ is a perfect square then $x\equiv0\bmod4$ or $x\equiv1\bmod4$~\cite[p.~21]{Aigner_Ziegler_2014}.
\item Cases 9, 10, 22 and 38 is unsolvable due to the assumed inequality $p<q$.
\item Cases 15 and 55 cannot be solved: Since $a\in\mathbb{Z}$, $b,p,q,\in\mathbb{N}$ and $(a,b)=1$, the conditions $q=pb^2-\frac{a^2}{b}$ and $pq=b^2-\frac{a^2}{b}$ require $b=1$ leading to $q=p-a^2$ and $pq=1-a^2$ which both have no solution.
\item Cases 16, 30, 31, 39, 46, 60 are impossible to solve: $p$ and $q$ are odd primes and the difference or sum of two even integers cannot be odd. To apply the argument for case 31, we first multiply its condition by $4b^4$.
\end{enumerate}

\section{Redundant conditions}
\label{app:reasons_redundancy}
We provide reasons for all conditions in Table~\ref{tab:cases_all} that are redundant: 
\begin{enumerate}[label=\alph*)]
\item Case 8: $\boldsymbol{p=qb^2-\frac{a^2}{b}}$: Due to $b,p,q$ are integers, the fraction $\frac{a^2}{b}$ needs to be an integer. Because $(a,b)=1$, solutions exist only if $b=1$. The condition simplifies to $p=q-a^2$, which is a special case of case 17. 
\item Case 14 and Case 44: $\boldsymbol{p=\frac{qb\pm2a^2}{4b^3}}$: 
We rewrite the conditions as $q=4b^2p\mp \frac{2a^2}{b}$. Due to $b,p,q$ are integers, the fraction $\frac{2a^2}{b}$ needs to be an integer. Because $(a,b)=1$, solutions exist only if $b=1$ or $b=2$. 
For case 14, setting $b=1$ leads to $q=4p-2a^2$ which is a special case of case 46 and is impossible to solve due to reason f in Appendix~\ref{app:reasons_unsatisfiability}. For case 14, setting $b=2$ leads to $q=2^4p-a^2$ which is a special case of case 47.
For case 44, setting $b=1$ leads to $q=4p+2a^2$ which is a special case of case 16 and is impossible to solve due to reason f in Appendix~\ref{app:reasons_unsatisfiability}. For case 14, setting $b=2$ leads to $q=2^4p+a^2$ which is a special case of case 17.
\item Cases 21 and 51: $\boldsymbol{q=\frac{pb\pm2a^2}{4b^3}}$: 
We rewrite the conditions as $p=4b^2q\mp \frac{2a^2}{b}$. Due to $b,p,q$ are integers, the fraction $\frac{2a^2}{b}$ needs to be an integer. Because $(a,b)=1$, solutions exist only if $b=1$ or $b=2$. 
For case 21, setting $b=1$ leads to $p=4q-2a^2$ which is a special case of case 39. For case 21, setting $b=2$ leads to $p=2^4q-a^2$ which is a special case of case 40. 
For case 51, setting $b=1$ leads to $p=4q+2a^2$ which is a special case of case 9. For case 21, setting $b=2$ leads to $p=2^4q+a^2$ which is a special case of case 10. Both cases are impossible to solve due to reason d in Appendix~\ref{app:reasons_unsatisfiability}.
\item Case 23: $\boldsymbol{pq=\frac{a^2+b}{b^3}}$: We rewrite the condition as $pqb^2=\frac{a^2}{b}+1$. Due to $b,p,q$ are integers, the fraction $\frac{a^2}{b}$ needs to be an integer. Because $(a,b)=1$, solutions exist only if $b=1$. Setting $b=1$ boils the condition down to $pq=a^2+1$ which is a special case of case 26.
\item Case 24: $\boldsymbol{pq=\frac{a^2+b^2}{b^2}}$: We rewrite the condition as $pq=\frac{a^2}{b^2}+1$. See the argumentation for case 23. Hence, setting $b=1$ boils the condition down to $pq=a^2+1$ which is a special case of case 26. 
\item Case 25: $\boldsymbol{pq=\frac{a^2+b^3}{b}}$: We rewrite the condition as $pq=\frac{a^2}{b}+b^2$. See the argumentation for case 23. Hence, setting $b=1$ boils the condition down to $pq=a^2+1$ which is a special case of case 26. 
\item Case 29 and Case 59: $\boldsymbol{pq=\frac{4b^3\pm2a^2}{b}}$: 
We rewrite the conditions as $pq=4b^2\pm\frac{2a^2}{b}$. 
Due to $b,p,q$ are integers, the fraction $\frac{2a^2}{b}$ needs to be an integer. 
Because $(a,b)=1$, solutions exist only if $b=1$ or $b=2$. For case 29, setting $b=1$ reduces the condition to $pq=2a^2+4$ which is a special case of case 33 that is impossible to solve due to reason b in Appendix~\ref{app:reasons_unsatisfiability}. For case 29, setting $b=2$ boils the condition down to $pq=a^2+2^4$ which is a special case of case 26. 
For case 59, setting $b=1$ reduces the condition to $pq=4-2a^2$ which is a special case of case 3 that is impossible to solve due to reason a in Appendix~\ref{app:reasons_unsatisfiability}. For case 59, setting $b=2$ boils the condition down to $pq=2^4-a^2$ which is a special case of case 56.
\item Case 45: $\boldsymbol{q=pb^2+\frac{a^2}{b}}$: 
Due to $b,p,q$ are integers, the fraction needs to be an integer. 
Because $(a,b)=1$, solutions exist only if $b=1$. This leads to $p=q-a^2$, which is a special case of case 17.
\item Case 52: $\boldsymbol{p=q-\frac{a^2}{b^2}}$: See argumentation for case 45.
\end{enumerate}
\newpage
\vspace{1em}
\bibliographystyle{unsrt}
\bibliography{main}

\end{document}